 \documentclass[final,1p]{elsarticle}
\usepackage{amssymb}
\biboptions{sort&compress}
\usepackage{amsthm}
\usepackage{bm}
\usepackage{amsmath}
\usepackage{enumerate}
\usepackage{enumitem}
 \usepackage{lineno}
\usepackage[caption=false]{subfig}
\newtheorem{theorem}{Theorem}[section]
\newtheorem{remark}[theorem]{Remark}
\newtheorem{lemma}[theorem]{Lemma}

\newtheorem{corollary}[theorem]{Corollary}
\newtheorem{proposition}[theorem]{Proposition}

\begin{document}

\begin{frontmatter}



\title{On the convergence of the nonlocal nonlinear model to the classical elasticity equation}


\author[label1]{H. A. Erbay\corref{cor1}}
    \ead{husnuata.erbay@ozyegin.edu.tr}
    \address[label1]{Department of Natural and Mathematical Sciences, Faculty of Engineering, Ozyegin University,  Cekmekoy 34794, Istanbul, Turkey}
    \cortext[cor1]{Corresponding author}

\author[label1]{S. Erbay}
    \ead{saadet.erbay@ozyegin.edu.tr}

\author[label2]{A. Erkip}
    \ead{albert@sabanciuniv.edu}
    \address[label2]{Faculty of Engineering and Natural Sciences, Sabanci University, Tuzla 34956, Istanbul,  Turkey}

\begin{abstract}
We consider  a general class of  convolution-type  nonlocal wave equations  modeling bidirectional propagation of nonlinear waves in a continuous medium.   In the limit of vanishing nonlocality we  study the behavior of solutions to the Cauchy problem.  We prove that, as the kernel functions of the convolution integral approach the Dirac delta function, the solutions  converge strongly to the corresponding solutions of the classical elasticity equation. An energy estimate with no loss of derivative plays a critical role in proving the convergence result. As a typical example, we consider the continuous limit of the discrete lattice dynamic model (the Fermi-Pasta-Ulam-Tsingou model) and  show that, as the lattice spacing approaches zero, solutions to the discrete lattice equation converge to the corresponding solutions of the classical elasticity equation.
\end{abstract}


\end{frontmatter}

%\linenumbers

\setcounter{equation}{0}
\section{Introduction}\label{sec1}

With a general class of kernel functions $\beta$,  the nonlocal equation
\begin{equation}
    u_{tt}=\beta \ast (u+g(u))_{xx},  \label{nw}
\end{equation}%
defines a class of  nonlinear wave equations modeling the bi-directional propagation of dispersive waves in a continuous media. Here $u=u(x,t)$ is a real-valued function,    $g$ is a sufficiently smooth nonlinear function satisfying $g(0)=g^{\prime }(0)=0$ and the symbol $\ast $ denotes convolution in the $x$-variable with the kernel function $\beta$. Throughout the manuscript we assume that the convolution  is a positive bounded operator on the Sobolev space $H^{s}(\mathbb{R})$ for any $s$. The linear dispersion relation $\xi \mapsto \omega^{2}(\xi)=\xi^{2}\widehat{\beta }(\xi)$,  where $\widehat{\beta }$ represents the Fourier transform of $\beta$, clearly demonstrates the dispersive nature of the solutions of (\ref{nw}).   For the particular case in which $\beta$ is the Dirac delta function, (\ref{nw}) becomes the classical elasticity equation
\begin{equation}
    u_{tt}=u_{xx}+g(u)_{xx}  \label{classical}
\end{equation}%
modeling the non-dispersive propagation of longitudinal waves in an elastic bar of infinite length, where a necessary and sufficient condition for it to be hyperbolic is $g^{\prime}(u)>-1$. The present work aims to establish  rigorously that, in the limit of vanishing nonlocality (that is,  as $\beta $ approaches the Dirac delta function),  the solution of (\ref{nw}) with a given initial data converges strongly  to the corresponding solution of (\ref{classical}) in a suitable Sobolev norm. More explicitly, we show that, for long waves of small-but-finite amplitude, the difference between the corresponding solutions of (\ref{nw}) and  (\ref{classical}) with the same initial data remains small over a long time scale.  As a concrete example  of the convergence result from the class (\ref{nw}) to (\ref{classical}), we consider the continuous limit of the Fermi–Pasta–Ulam-Tsingou (FPUT)   lattice  model describing the dynamics of the one-dimensional doubly-infinite chain of identical masses with nearest neighbor interactions. The FPUT  lattice dynamics model is also a member of the class (\ref{nw}), which corresponds to the triangular kernel to be defined later. We show that, as the lattice spacing approaches zero, the leading approximation of the discrete lattice model approaches the continuous model (\ref{classical}).

The nonlocal nonlinear wave equation (\ref{nw}) was introduced in \cite{Duruk2010} to model a one-dimen\-sional motion of a nonlocally and nonlinearly elastic medium. The dependent variable $u(x,t)$ is the strain (the spatial derivative of displacement) at the point $x$ at time $t$, the given function $g(u)$ is the nonlinear part of the (local) stress and (\ref{nw}) is the dimensionless equation of motion for the strain (i.e. Newton's second law of motion for constant mass density and in the absence of body forces). The theory of nonlocal elasticity \cite{Eringen2002} is an active and interesting subject of research in solid mechanics. In recent years, it has been widely used to predict a more realistic picture of material behavior at the micro/nano length scales. Its main difference with classical elasticity relies  on the  assumption that the stress  at a given  point is a function of the strain field at every point in the body, that is,  the nonlocal stress is assumed to be  the spatial average of the local stress field  over the body weighted with a kernel. This fact is reflected in (\ref{nw}) through the convolution integral term  which introduces  a length scale into the model. So  this  length scale  is responsible for dispersion of  wave solutions to (\ref{nw}).   It is worth mentioning that the class (\ref{nw}) of nonlocal nonlinear wave equations covers a variety of well-known dispersive wave equations, for instance  the improved Boussinesq equation $u_{tt}-u_{xx}-u_{xxtt}=g(u)_{xx}$ corresponding to the exponential kernel $\beta(x)=\frac{1}{2}e^{-\vert x \vert}$ is the  most widely known member of the class.   Depending on the weakness of dispersion, the class (\ref{nw}) provides a hierarchy of nonlinear dispersive wave equations which serve as approximations of (\ref{classical}).  It is natural to ask whether, if the dispersive effect is sufficiently weak,  solutions of (\ref{classical}) can be approximated using appropriate solutions of (\ref{nw}). So, in this paper, we concentrate on proving  the convergence of the nonlocal nonlinear model (\ref{nw}) to the classical elasticity equation (\ref{classical}) on a long time interval in the zero-dispersion limit.

A collection of articles \cite{Erbay2011a,Duruk2011,Erbay2015,Erbay2018,Erbay2019,Erbay2020}  focuses on  various aspects of mathematical analysis of the initial-value problems associated with (\ref{nw}). In those studies it is assumed that  the Fourier transform $\widehat{\beta }$ of the kernel function $\beta$ satisfies the ellipticity and boundedness condition
\begin{equation}
    0\leq \widehat{\beta }(\xi )\leq C(1+\xi ^{2})^{-\frac{r}{2}}\text{ \ for some } r\geq 0.  \label{r}
\end{equation}%
Here we remark that the necessary condition for wave propagation is given by $\widehat{\beta }(\xi )\geq 0$ and that general dispersive properties of waves are characterized by the kernel function. The exponent $r$ in (\ref{r})  denotes the rate of decay of $\widehat{\beta }$ as $\vert \xi \vert \rightarrow \infty$  and it is determined by  the smoothness of $\beta $ (faster decay reflects a greater smoothness).  The limiting case $r=0$ corresponds to the case when $\beta$ is integrable or more generally a finite measure on $\mathbb{R}$. We note that $\widehat{\beta }$ is real and hence $\beta$ should be an even function.

The local well-posedness of the Cauchy problem for (\ref{nw}) was proved in \cite{Duruk2010} under the strong smoothness assumption $r\geq 2$. This constraint is equivalent to say that  the operator $\beta^{\prime\prime}\ast (\cdot)$ to be bounded on $H^{s}(\mathbb{R})$. As a result,  (\ref{nw}) becomes a Banach space valued ODE.  In \cite{Erbay2019}  the well-posedness result was extended to  nonsmooth kernels satisfying  $0\leq r <2$ for which (\ref{nw}) is genuinely hyperbolic; hence requires a totally different approach from that in the $r\geq 2$ case. Another improvement in \cite{Erbay2019} was a long-time existence result in the small amplitude regime.   The main  tool used there was a Nash-Moser iteration scheme necessitated by the loss of derivative in the energy estimate.  In this work we first give an alternative proof for the local well-posedness of the Cauchy problem for (\ref{nw}) with $r\geq 0$. The  improvement here is   due to the construction of a quasilinear energy functional allowing  a better energy estimate which does not involve a loss of derivative.  This allows us to use a Picard  iteration scheme to prove the local well-posedness result for (\ref{nw}). This result is also stronger than the one in  \cite{Erbay2019} as it does not require extra smoothness of the initial data, which was due to the Nash-Moser scheme there. Furthermore, this  improved energy estimate   allows us to establish the blow-up criteria for nonsmooth kernels, extending of the blow-up result obtained in \cite{Duruk2010} smooth kernels (i.e. for the ODE case). We show that the maximal existence time is controlled  by the $L^{\infty }$  norms of $u$, $u_{x}$ and $u_{t}$ (that may lead to wave breaking) in contrast to the $L^{\infty }$  norm of $u$ in the $r\geq 2$ case.

In the present work we prove that the  limit of (\ref{nw}) is  (\ref{classical}) in the sense that  solutions of (\ref{nw}) approximate to the solution of (\ref{classical})  in the limit of vanishing nonlocality. This is inspired by the comparison result established in  \cite{Erbay2020} where two solutions  corresponding to two different kernels were compared  in the long-wave limit. Here we concentrate on the particular case where one of the kernels is the Dirac delta function corresponding to the zero-dispersion case. In other words we compare  (\ref{nw}) for  an arbitrary kernel with  the classical elasticity equation. The main step is to obtain  uniform estimates of solutions with respect to the parameter characterizing dispersion and to estimate the difference between the solutions  using the energy estimate. Nevertheless, the fact that our energy estimate does not involve loss of derivative allows us to get a stronger convergence result than that in \cite{Erbay2020}.

In this work we will consider the kernels that satisfy
\begin{equation}
    0\leq \widehat{\beta }(\xi )\leq C  \label{rbound}
\end{equation}%
with $C=\Vert \widehat{\beta }\Vert_{L^{\infty}}$. We note that any kernel satisfying (\ref{r}) with some $r\geq 0$ also satisfies the weaker estimate (\ref{rbound}). In other words, we consider all integrable kernels as well as all finite measures on  $\mathbb{R}$. This is the largest class of kernels that we will consider  throughout this work. Clearly, one can obtain sharper results assuming extra smoothness on $\beta$ (i.e. $r>0$). We will indicate those cases in our local existence theorems.

The plan of this paper is as follows. In Section \ref{sec2} we write (\ref{nw}) as a first-order nonlinear system and prove an energy estimate for strong solutions of the corresponding linearized system. In Sections \ref{sec3} and \ref{sec4} we prove the  well-posedness of the Cauchy problem for the linearized system and the nonlinear system, respectively.  Section \ref{sec5} is concerned with the dependence of the existence time of the solutions on the initial data and the size of the nonlinear term, which makes possible to  extend the local well-posedness result to long time intervals. In the same section, as a by-product, we prove  the result stating that the finite-time blow-up of solutions can occur only in the $L^{\infty }$ norms of  $u$, $u_{x}$ and $u_{t}$. In Section \ref{sec6},  we are able to prove the main result of the paper stating that, in limit of vanishing nonlocality (i.e. in the zero-dispersion limit),  small-but-finite amplitude solutions of the Cauchy problem for (\ref{nw}) converge to  the corresponding solution of (\ref{classical}). In Section \ref{sec7} we provide a rigorous justification of the convergence of the FPUT discrete lattice dynamic equation corresponding to the triangular kernel  to the classical elasticity equation as the lattice spacing approaches zero.

Throughout the paper, we use the standard notation for function spaces. The Fourier transform $\widehat{u}$ of $u$ is defined by $\widehat{u}(\xi
)=\int_{\mathbb{R}}u(x)e^{-i\xi x}dx$. The norm of $u$ in the Lebesgue space $L^{p}(\mathbb{R})$ ($1\leq p\leq \infty $) is represented by $\Vert \cdot\Vert_{L^{p}}$. To denote the inner product of $u$ and $v$ in $L^{2}$ space, the symbol $\langle u,v\rangle $ is used. The notation $H^{s}=H^{s}(\mathbb{R})$ (for $s\in\mathbb{R}$) is used to denote the $L^{2}$ -based Sobolev space of order $s$ on $\mathbb{R}$, with the norm $\Vert u\Vert _{H^{s}}=\big( \int_{\mathbb{R}}(1+\xi ^{2})^{s}|\widehat{u}(\xi )|^{2}d\xi \big)^{1/2}$.  $C$ is a generic positive constant. Partial differentiations are denoted by $D_{x}$ etc. Further, to simplify the presentation in the rest of the paper,  we introduce the notations $X^{s}$ and $Y^{s_{1},s_{2}}$ to refer the spaces  defined by
\begin{equation}
    X^{s}=C\big([0, T], H^{s}\big),~~~~Y^{s_{1},s_{2}}=C\big([0, T], H^{s_{1}}\big) \cap C^{1}\big([0, T], H^{s_{2}}\big), \label{XYuzay}
\end{equation}
for fixed $T>0$.  Here the spaces $X^{s}$ and $Y^{s_{1},s_{2}}$ are endowed with the usual norms
\begin{equation}
    \Vert u\Vert_{X^{s}}=\sup\limits_{0\leq t\leq T}\Vert u(t)\Vert _{H^{s}}, ~~~~
    \Vert u\Vert _{Y^{s_{1},s_{2}}}
        =\sup\limits_{0\leq t\leq T}\Vert u(t)\Vert_{H^{s_{1}}}+\sup\limits_{0\leq t\leq T}\Vert u_{t}(t)\Vert_{H^{s_{2}}}, \label{XYnorm}
\end{equation}
respectively.

\setcounter{equation}{0}
\section{An energy estimate}\label{sec2}

The proof of the local well-posedness of (\ref{nw}) in Sobolev spaces $H^{s}$ is based on the energy estimates and Gronwall's inequality. In this section, we construct an energy estimate for the first-order inhomogeneous linearized system related to (\ref{nw}). Throughout this section we assume that $\beta$ satisfies (\ref{rbound}).

In the rest of the work, we will adopt the notations $\Lambda^{s}=(1-D_{x}^{2})^{s/2}$ and $[\Lambda^{s},f] g=\Lambda^{s}(f g)-f\Lambda^{s}g$. Furthermore, to prove our estimates below and in the next sections, we will need the following commutator estimates \cite{Kato1988}:
\begin{lemma}\label{lem2.1}
    Let  $s>0$. Then for all $f, g$ satisfying $f\in H^{s}$, $D_{x}f\in L^{\infty}$, $g\in H^{s-1}\cap L^{\infty}$,
    \begin{equation*}
        \big\Vert [ \Lambda^{s},f] g\big\Vert_{L^{2}}
        \leq C \big(\Vert D_{x}f\Vert_{L^{\infty}}\Vert g\Vert_{H^{s-1}}+\Vert f\Vert _{H^{s}}\Vert g\Vert_{L^{\infty }}\big).
    \end{equation*}%
    In particular, when $s>3/2$, due to the Sobolev embeddings $H^{s-1}\subset L^{\infty}$, for all $f, g\in H^{s}$
    \begin{equation*}
        \big\Vert [ \Lambda^{s},f] D_{x}g\big\Vert_{L^{2}}\leq C \Vert f\Vert_{H^{s}}\Vert g\Vert_{H^{s}}.
    \end{equation*}%
\end{lemma}

We first define the  pseudo-differential operator $K$  by
\begin{equation}
    Ku=\mathcal{F}^{-1}\big( k(\xi )\widehat{u}(\xi)\big) \text{ with \ } k(\xi )=\big(\widehat{\beta}(\xi)\big)^{1/2}  \label{K1}
\end{equation}
for any $\beta$ satisfying (\ref{rbound}). Clearly, $K$ is bounded on $H^{s}$ with operator norm $C=\Vert \widehat{\beta }\Vert_{L^{\infty}}^{1/2}$.   Then  the Cauchy problem
\begin{eqnarray}
    && u_{t}=Kv_{x},\text{ \ \ \ \ \ \ \ \ \ \ \ \ \ \ \ \ \ \ \ }u(x,0)=u_{0}(x),~~\label{systema} \\
    && v_{t}=K\big((1+g^{\prime }(u))u_{x}\big),\text{ \ \ \ }v(x,0)=v_{0}(x)  \label{systemb}
\end{eqnarray}%
is equivalent to the Cauchy problem defined by (\ref{nw})  and the initial data $u(x,0)=u_{0}(x)$,  $u_{t}(x,0)=K(v_{0}(x))_{x}$.  We now study the corresponding  linearized problem
\begin{eqnarray}
    && u_{t}=Kv_{x},\text{ \ \ \ \ \ \ \ \ \ \ \ \ \ \ \ \ \ \ }u(x,0)=u_{0}(x),~~ \label{linsysa} \\
    && v_{t}=K\big((1+w)u_{x}\big),\text{ \ \ \ \ \ \ }v(x,0)=v_{0}(x).  \label{linsysb}
\end{eqnarray}
Here  $w$ is a  given fixed function and the hyperbolicity of the linearized system (\ref{linsysa})-(\ref{linsysb})  is guaranteed by the condition $w(x,t)>-1$. We define the $H^{s}$-energy functional of (\ref{linsysa})-(\ref{linsysb}) by
\begin{equation}
    \mathcal{E}_{s}^{2}(t)
        =\frac{1}{2}\int_{\mathbb{R}} \Big( ( 1+w(x,t))(\Lambda^{s}u(x,t))^{2}+(\Lambda^{s}v(x,t))^{2}\Big) dx,  \label{energy}
\end{equation}%
where $\Lambda^{s}=(1-D_{x}^{2})^{s/2}$.  Note that when $w$ is assumed to satisfy
\begin{equation}
    0<c_{1}\leq 1+w(x, t)\leq c_{2}  \label{hyperbolic}
\end{equation}%
where $c_{1}$ and $c_{2}$ are constants, $\mathcal{E}_{s}^{2}(t)$ is equivalent to the norm $\Vert u(t)\Vert_{H^{s}}^{2}+\Vert v(t)\Vert _{H^{s}}^{2}$.

\begin{lemma}\label{lem2.2}
    Let $s>3/2$, $X^{s}=C\big([0, T], H^{s}\big)$, $Y^{s, s-1}=C\big([0, T], H^{s}\big) \cap C^{1}\big([0, T], H^{s-1}\big)$. Let $u_{0},v_{0}\in H^{s}$ and $w\in Y^{s,s-1}$ with $~0<c_{1}\leq 1+w(x,t)~$    for all $( x,t) \in \mathbb{R}\times [0, T] $. Suppose $u, v\in X^{s}$ satisfy the inhomogeneous system
    \begin{eqnarray}
        && u_{t}=Kv_{x}~+F_{1},\text{ \ \ \ \ \ \ \ \ \ \ \ \ \ \ \ \ \ \ \ \ \  }    u(x,0)=u_{0}(x),  \label{inhomA} \\
        && v_{t}=K\big((1+w)u_{x}\big)+F_{2}+G,\text{ \ \ \ \ \ }v(x,0)=v_{0}(x)  \label{inhomB}
    \end{eqnarray}%
    on $\mathbb{R}\times [ 0,T]$, where $G$ satisfies $\Vert G(t)\Vert_{H^{s}}\leq C \Vert u(t)\Vert _{H^{s}}$ and $F_{1},F_{2}\in X^{s}$. Then we have the estimate
    \begin{equation}
         \Vert u(t)\Vert_{H^{s}}+\Vert v(t)\Vert_{H^{s}}
           \leq  \big(\Vert u_{0}\Vert_{H^{s}}+\Vert v_{0}\Vert_{H^{s}}\big)e^{Ct}
            +C(\Vert F_{1}\Vert _{X^{s}}+\Vert F_{2}\Vert_{X^{s}})(e^{Ct}-1)  \label{energy21}
    \end{equation}%
    for $0\leq t\leq T$,     where $C$ is a constant dependent on $\Vert w\Vert _{Y^{s,s-1}}$.
\end{lemma}
\begin{proof}
    Throughout the proof we will use the abbreviation $\widetilde{w}=1+w$ to simplify the presentation.     Differentiating (\ref{energy})  and using (\ref{inhomA})-(\ref{inhomB}) we get
    \begin{eqnarray}\label{xx}
        \frac{d}{dt}\mathcal{E}_{s}^{2}(t)
            &=&\int_\mathbb{R} \Big (\frac{1}{2} w_{t}(\Lambda^{s}u) ^{2}
                +\widetilde{w}(\Lambda^{s}u_{t})(\Lambda^{s}u) +(\Lambda^{s}v_{t})(\Lambda^{s}v) \Big) dx \notag \\
            &=&\int_\mathbb{R} \Big(\frac{1}{2}w_{t}( \Lambda ^{s}u)^{2}+\widetilde{w}(\Lambda ^{s} F_{1})(\Lambda ^{s}u)
                +(\Lambda ^{s}F_{2})(\Lambda^{s}v)+(\Lambda ^{s}G)(\Lambda^{s}v)   \notag \\
            &&  ~~~~~+\widetilde{w} (\Lambda ^{s}K v_{x})(\Lambda ^{s}u) +(\Lambda^{s}K(\widetilde{w} u_{x}))(\Lambda^{s}v) \Big)dx.
    \end{eqnarray}
    We handle the last two terms on the right-hand-side of (\ref{xx}) separately:
    \begin{eqnarray}
         I(t) &=& \big\langle \widetilde{w} (\Lambda ^{s}K v_{x}), \Lambda^{s}u\big\rangle
            +\big\langle\Lambda ^{s}K(\widetilde{w} u_{x}),\Lambda ^{s}v\big\rangle\notag\\
           &=& \big\langle D_{x}(\widetilde{w} \Lambda ^{s}K v)-\widetilde{w}_x(\Lambda ^{s}K v), \Lambda ^{s}u\big\rangle
            +\big\langle\Lambda ^{s}K(\widetilde{w} u_{x}),\Lambda ^{s}v\big\rangle\notag\\
           &=& -\big\langle \widetilde{w} (\Lambda ^{s}K v), \Lambda ^{s}u_x\big\rangle
                -\big\langle \widetilde{w}_x (\Lambda ^{s}K v), \Lambda ^{s}u\big\rangle
                +\big\langle\Lambda ^{s}(\widetilde{w} u_{x}),\Lambda ^{s}K v\big\rangle\notag\\
           &=& -\big\langle w_x (\Lambda ^{s}K v), \Lambda^{s}u\big\rangle
                +\big\langle \Lambda^{s}K v, -\widetilde{w}(\Lambda ^{s}u_x)+\Lambda^{s}(\widetilde{w} u_{x}) \big\rangle \notag\\
           &=& -\big\langle w_x (\Lambda ^{s}K v), \Lambda^{s}u\big\rangle
                +\big\langle \Lambda ^{s}K v, [\Lambda ^{s},w] u_{x} \big\rangle,
    \end{eqnarray}
     with $[\Lambda^{s},\widetilde{w}] u_{x}=[\Lambda^{s},w] u_{x}=\Lambda^{s}(w u_{x})-w\Lambda^{s}u_{x}$,  where we have used the facts that the operators $\Lambda^{s}$, $K$ and $D_{x}$ all commute and $K$ is self-adjoint on $ L^{2}$. Thus, (\ref{xx}) takes the form
    \begin{eqnarray*}
        \frac{d}{dt}\mathcal{E}_{s}^{2}(t)
            &=&\int_\mathbb{R} \Big( \frac{1}{2}w_{t}( \Lambda ^{s}u) ^{2}- \widetilde{w}_x (\Lambda ^{s}K v)(\Lambda ^{s}u)
                +(\Lambda^{s}K v) [\Lambda^{s},w] u_{x}  \\
            &~& ~~~~~+\widetilde{w}(\Lambda ^{s} F_{1})(\Lambda ^{s}u)+(\Lambda ^{s}F_{2})(\Lambda^{s}v)+(\Lambda ^{s}G)(\Lambda^{s}v)\Big)dx.
    \end{eqnarray*}
    Since $\Vert Ku \Vert_{H^{s}}\leq C\Vert u \Vert_{H^{s-\frac{r}{2}}}\leq C\Vert u \Vert_{H^{s}}$, we deduce from the above equation that
    \begin{eqnarray}
        && \!\!\!\!\!\!\!\!\!\!\!
      \frac{d}{dt}\mathcal{E}_{s}^{2}
            \leq  C\Big(\Vert w_{t}\Vert_{L^{\infty }} \Vert u\Vert _{H^{s}}^{2}
                +\Vert w_{x}\Vert_{L^{\infty}} \Vert v\Vert_{H^{s}}\Vert u\Vert_{H^{s}}
                +\Vert v\Vert_{H^{s}} \big\Vert [ \Lambda ^{s},w] u_{x}\big\Vert_{L^{2}}
                 \notag \\
            &&~~~~~~~~  +\Vert \widetilde{w}\Vert_{L^{\infty}}\Vert F_{1}\Vert _{H^{s}}\Vert u\Vert_{H^{s}}
                +\Vert F_{2}\Vert _{H^{s}}\Vert v\Vert _{H^{s}}
                +\Vert G\Vert _{H^{s}}\Vert v\Vert _{H^{s}}\Big).  \label{e-com}
    \end{eqnarray}
    Using  the relation $\Vert G(t)\Vert_{H^{s}}\leq C \Vert u(t)\Vert _{H^{s}}$ and the commutator estimate of Lemma \ref{lem2.1}
    \begin{equation*}
        \big\Vert [ \Lambda^{s},w] u_{x}\big\Vert _{L^{2}}
            \leq  C\Vert w\Vert _{H^{s}}\Vert u\Vert _{H^{s}},
    \end{equation*}
    we get
    \begin{eqnarray}
        && \!\!\!\!\!\!\!\!\!\!\!
      \frac{d}{dt}\mathcal{E}_{s}^{2}
            \leq  C\Big(\Vert w_{t}\Vert _{L^{\infty }} \Vert u\Vert _{H^{s}}^{2}
                +\Vert w_{x}\Vert _{L^{\infty}} \Vert v\Vert _{H^{s}}\Vert u\Vert _{H^{s}}
                +\big(1+\Vert w\Vert_{H^{s}}\big)\Vert v\Vert_{H^{s}} \Vert u\Vert_{H^{s}}
                 \notag \\
            &&~~~~~~~~~~~~  + \Vert \widetilde{w}\Vert_{L^{\infty}} \Vert F_{1}\Vert _{H^{s}}\Vert u\Vert_{H^{s}}
                +\Vert F_{2}\Vert _{H^{s}}\Vert v\Vert _{H^{s}} \Big).     \label{energy-estimate}
    \end{eqnarray}
    For $s>3/2$, we have
    \begin{equation*}
        \sup_{0\leq t\leq T}\Big(\Vert w_{t}(t)\Vert _{L^{\infty}}+\Vert w_{x}(t)\Vert _{L^{\infty }}+\Vert w(t)\Vert_{H^{s}}+1\Big)
        \leq C_{1},
    \end{equation*}%
    and
    \begin{equation*}
        \sup_{0\leq t\leq T}\Vert \widetilde{w}(t)\Vert _{L^{\infty}}    \leq C_{2}
    \end{equation*}%
    where $C_{1}$ and $C_{2}$ are constants depending on $\Vert w\Vert _{Y^{s,s-1}}$. So, using the fact that $\mathcal{E}_{s}(t)$ is equivalent to the norm $\Vert u(t)\Vert_{H^{s}}+\Vert v(t)\Vert_{H^{s}}$,     (\ref{energy-estimate}) becomes
    \begin{equation*}
        \frac{d}{dt}\mathcal{E}_{s}^{2}(t)
            \leq C_{1}\mathcal{E}_{s}^{2}(t)+C_{2}\big(\Vert F_{1}\Vert _{H^{s}}+\Vert F_{2}\Vert_{H^{s}}\big)\mathcal{E}_{s}(t).
    \end{equation*}%
    By Gronwall's inequality we get,
    \begin{equation*}
        \mathcal{E}_{s}(t)
        \leq \mathcal{E}_{s}(0)e^{C_{1}t}+C_{2}\big(\Vert F_{1}\Vert_{X^{s}}+\Vert F_{2}\Vert _{X^{s}}\big)\big(e^{C_{1}t}-1\big).
    \end{equation*}%
    This completes the proof of (\ref{energy21}).
\end{proof}

\setcounter{equation}{0}
\section{Well-posedness for the linear system}\label{sec3}

In this section we establish the  well-posedness of  the Cauchy problem (\ref{linsysa})-(\ref{linsysb}) defined for the linearized system.
\begin{lemma}\label{lem3.1}
        Suppose $\beta$ satisfies (\ref{r}) for some $r\geq 0$.     Let $s>5/2$, $\sigma=\min \{ s,s+\frac{r}{2}-1\}$, $X^{s}=C\big([0, T], H^{s}\big)$ and $Y^{s, \sigma}=  C\big( [0, T] ,H^{s}\big) \cap C^{1}\big( [0,T],H^{\sigma}\big)$.  Let $u_{0},v_{0}\in H^{s}$ and $w\in Y^{s,\sigma}$ with $~0<c_{1}\leq 1+w(x,t)~$
    for all $(x,t) \in \mathbb{R}\times [0, T] $. Then there exist unique $u,v\in Y^{s,\sigma}$ satisfying (\ref{linsysa})-(\ref{linsysb}) on $\mathbb{R}\times [ 0,T]$.
\end{lemma}
\begin{proof}
    For the existence proof of the linearized system, we follow the standard hyperbolic approach \cite{Taylor2011} in which a solution of (\ref{linsysa})-(\ref{linsysb}) is obtained as a limit of solutions to
    \begin{eqnarray}
        u_{t}^{h} &=&J^{h}Kv_{x}^{h},~\text{\ \ \ \ \ \ \ \ \ \ \ \ \ \ }u^{h}(x,0)=u_{0}(x),  \label{molA} \\
        v_{t}^{h} &=&J^{h}K\big((1+w)u_{x}^{h}\big),\text{ \ \ }v^{h}(x,0)=v_{0}(x).  \label{molB}
    \end{eqnarray}%
    Here $J^{h}$ is a Friederichs mollifier given by
    \begin{equation*}
        J^{h}\varphi(x)=\frac{1}{h}\int_{\mathbb{R}} \eta (\frac{x-y}{h})\varphi(y)dy
    \end{equation*}
    with some nonnegative $\eta \in C_{0}^{\infty }(\mathbb{R)}$ and $\int_{\mathbb{R}} \eta(x) dx=1.$ Due to the regularizing effect of $J^{h},$ (\ref{molA})-(\ref{molB}) is an $X^{s}$-valued ODE system and hence has solution $u^{h},v^{h}\in X^{s}$. Moreover, since $J^{h}K$ satisfies the same bounds as $K$, the energy estimate of Lemma \ref{lem2.2} will hold uniformly in $h$. It follows that $(u^{h},v^{h})$ is uniformly bounded in $X^{s}$. We will show that for any sequence $h_{n}\rightarrow 0$, the solutions $(u^{h_{n}},v^{h_{n}})$ form Cauchy sequences in a lower norm. The key ingredient in the proof is the mollifier estimate \cite{Mats2017}
    \begin{equation}
        \Vert J^{h_{n}}\varphi -J^{h_{m}}\varphi \Vert _{H^{s-1}}
        \leq \vert h_{n}-h_{m}\vert \,\Vert \varphi \Vert _{H^{s}}. \label{molest}
    \end{equation}%
   On the other hand we have the direct estimate $\Vert J^{h_{n}}\varphi -J^{h_{m}}\varphi \Vert_{H^{s}}
        \leq C \Vert \varphi \Vert _{H^{s}}$. Using the Gagliardo–Nirenberg type inequality \cite{Brezeis2018} for any $0 <\alpha <1$ we have
   \begin{equation}
        \Vert J^{h_{n}}\varphi -J^{h_{m}}\varphi \Vert _{H^{s-\alpha}}
        \leq C \vert h_{n}-h_{m}\vert^{\alpha} \Vert \varphi \Vert _{H^{s}}. \label{molestek}
    \end{equation}%
    Given $s>5/2$ we first choose $0 <\alpha <1$ so that $s-\alpha>5/2$. We start by writing (\ref{molA})-(\ref{molB})  for $h_{n}$ and $h_{m}$ and setting $p=u^{h_{n}}-u^{h_{m}}$, $q=v^{h_{n}}-v^{h_{m}}$. Then, (\ref{molA})-(\ref{molB}) imply that the differences $p,q$ satisfy the initial value problem
    \begin{eqnarray*}
        p_{t} &=&J^{h_{n}}Kq_{x}+(J^{h_{n}}-J^{h_{m}})Kv_{x}^{h_{m}},~\text{\ }
        ~\text{\ \ \ \ \ \ \ \ \ \ \ \ \ \ \ \ \ \ \ \ \ \ \ \ \ \ \ }p(x,0)=0,   \\
        q_{t} &=&J^{h_{n}}K\big((1+w)p_{x}\big)+(J^{h_{n}}-J^{h_{m}})K\big((1+w)u_{x}^{h_{m}}\big),~\text{\ }~\text{ \ }q(x,0)=0.
    \end{eqnarray*}
    The above problem is in the form of (\ref{inhomA})-(\ref{inhomB}) if we take
    \begin{equation*}
    F_{1}=(J^{h_{n}}-J^{h_{m}})Kv_{x}^{h_{m}}, ~~~~F_{2}=(J^{h_{n}}-J^{h_{m}})K\big((1+w)u_{x}^{h_{m}}\big), ~~~~G=0.
    \end{equation*}
    By (\ref{molestek})  we have
    \begin{eqnarray*}
       && \Vert F_{1}\Vert_{X^{s-1-\alpha}}\leq C\vert h_{n}-h_{m}\vert^{\alpha} \Vert v\Vert_{H^{s}} \leq C\vert h_{n}-h_{m}\vert^{\alpha}, \\
       && \Vert F_{2}\Vert_{X^{s-1-\alpha}}\leq C\vert h_{n}-h_{m}\vert^{\alpha} \Vert u\Vert_{H^{s}} \leq C\vert h_{n}-h_{m}\vert^{\alpha} .
    \end{eqnarray*}
      Then the energy estimate (\ref{energy21}) of  Lemma \ref{lem2.2} implies that the solutions $(p,q)=\big(u^{h_{n}}-u^{h_{m}}, v^{h_{n}}-v^{h_{m}}\big)$ satisfy
    \begin{equation*}
        \Vert u^{h_{n}}(t)-u^{h_{m}}(t)\Vert_{H^{s-1-\alpha}}+\Vert v^{h_{n}}(t)-v^{h_{m}}(t)\Vert _{H^{s-1-\alpha}}
        \leq C\vert h_{n}-h_{m}\vert^{\alpha} (e^{Ct}-1).
    \end{equation*}%
    Hence $\lim_{h \rightarrow 0}( u^{h},v^{h}) =(u,v)$ exists in $X^{s-1-\alpha}$.  On the other hand, $(u^{h},v^{h}) $ being bounded in $X^{s}$ will have a weak limit in $X^{s}$.  Uniqueness of the limit shows regularity, namely that $u,v\in X^{s}$ indeed satisfies (\ref{linsysa})-(\ref{linsysb}). The second part of the assertion, $u_{t},v_{t}\in C\big([0, T], H^{\sigma}\big)$ (with $\sigma=s+\frac{r}{2}-1$ when $r<2$), follows directly from the smoothing effect of $K$  in (\ref{molA})-(\ref{molB}). Finally, uniqueness for the linear system follows directly from the energy estimate.
\end{proof}

\setcounter{equation}{0}
\section{Local well-posedness for the nonlinear system}\label{sec4}

Once having proved the  well-posedness of the linearized system,  the next stage is to prove the local well-posedness of  the nonlinear system  (\ref{systema})-(\ref{systemb}). In this section we use Picard's iterations to prove the local well-posedness of (\ref{systema})-(\ref{systemb}). In the process, we will make use of the nonlinear estimates (see \cite{Alinhac2007, Constantin2002}):
\begin{lemma}\label{lem4.1}
    Let $h\in C^{\infty }(\mathbb{R})$ with $h(0)=0$. Then, for any $s\geq 0$ and $u,v\in L^{\infty}\cap H^{s}$,
    \begin{enumerate}
    \item $h(u)\in H^{s}$ with $\ \Vert h(u)\Vert _{H^{s}}\leq C_{1}\Vert u\Vert _{H^{s}}$ where $C_{1}$ depends on $h$ and $\Vert u\Vert _{L^{\infty }}.$
    \item $\Vert h(u)-h(v)\Vert _{H^{s}}\leq C_{2}\Vert u-v\Vert _{H^{s}}$ where $C_{2}$ depends on $h$ and $\Vert u\Vert_{L^{\infty }},\Vert v\Vert_{L^{\infty}},\Vert u\Vert_{H^{s}},$ and $\Vert v\Vert _{H^{s}}.$
    \end{enumerate}
\end{lemma}

To employ the energy estimate of Lemma \ref{lem2.2}, we want the initial value $u(x,0)=u_{0}(x)$ to satisfy
\begin{equation*}
    0<c_{1}\leq 1+g^{\prime }(u_{0}(x))
\end{equation*}%
for all $x\in \mathbb{R}$. This can be achieved as follows: Since $g^{\prime}(0)=0,$ one has $0<c_{1}\leq 1+g^{\prime }(z)$ for sufficiently small $\vert z\vert .$ By the Sobolev embedding theorem, this in turn implies that there is some $\gamma $ so that $0<c_{1}\leq 1+g^{\prime
}(z(x,t))$ whenever $\Vert z(t)\Vert _{H^{s}}\leq \gamma$.

We start by assuming that $\Vert u_{0}\Vert _{H^{s}}+\Vert v_{0}\Vert
_{H^{s}}\leq \frac{\gamma }{2}$ so that \ $1+g^{\prime }(u_{0})\geq c_{1}$. We inductively define the iterates  $(u^{n+1},v^{n+1})$ as the solution of the linear system
\begin{eqnarray*}
    u_{t}^{n+1} &=&Kv_{x}^{n+1},~\text{ \ \ \ \ \ \ \ \ \ \ \ \ \ \ \ \ \ \ }u^{n+1}(x,0)=u_{0}(x), \\
    v_{t}^{n+1} &=&K\big((1+g^{\prime }(u^{n}))u_{x}^{n+1}\big),\text{ \ }v^{n+1}(x,0)=v_{0}(x)
\end{eqnarray*}%
with $(u^{0},v^{0})=(u_{0},v_{0})$. By the energy estimate (\ref{energy21}) of Lemma \ref{lem2.2} for $F_{1}=F_{2}=G=0$  and $t\leq T_{0}=\frac{\log 2}{C}$ we have
\begin{equation*}
    \Vert u^{1}(t)\Vert _{H^{s}}+\Vert v^{1}(t)\Vert _{H^{s}}
        \leq \big(\Vert u_{0}\Vert _{H^{s}}+\Vert v_{0}\Vert _{H^{s}}\big)e^{Ct}\leq \frac{\gamma}{2}e^{Ct}\leq \gamma
\end{equation*}%
with $C=C\big(\Vert g^{\prime }(u_{0})\Vert _{H^{s}}\big)\leq C(\gamma)$. A similar estimate will also work for all $(u^{n},v^{n}),$ hence we have $0<c_{1}\leq 1+g^{\prime }(u^{n}(x,t))$ \ for all $x$ and $t\leq T_{0}=\frac{\log 2}{C}$. We next estimate $( p^{n+1}, q^{n+1}) =(
u^{n+1}-u^{n},v^{n+1}-v^{n})$. Then, the differences $p^{n+1}, q^{n+1}$ satisfy
\begin{eqnarray*}
    p_{t}^{n+1} &=&Kq_{x}^{n+1},
    ~\text{\ \ \ \ \ \ \ \ \ \ \ \ \ \ \ \ \ \ \ \ \ \ \ \ \ \ \ \ \ \ \ \ \ \ \ \ \ \ \ \ \ \ \ \ \ \ \ \ \ \ \ \ \ \ \ \ \ \ \ \ }p^{n+1}(x,0)=0, \\
    q_{t}^{n+1} &=&K\big((1+g^{\prime }(u^{n}))p_{x}^{n+1}\big)+K\big((g^{\prime}(u^{n})-g^{\prime }(u^{n-1}))u_{x}^{n}\big),~\text{\ }~\text{ \ }q^{n+1}(x,0)=0.
\end{eqnarray*}%
Applying Lemma \ref{lem2.2} for $s-1$ rather than $s$, with $F_{1}=G=0$ and $F_{2}=K\big(g^{\prime}(u^{n})-g^{\prime }(u^{n-1})\big)u_{x}^{n}$, we get
\begin{eqnarray*}
    \Vert p^{n+1}(t)\Vert _{H^{s-1}}+\Vert q^{n+1}(t)\Vert_{H^{s-1}}
        &\leq & C \big\Vert K\big(g^{\prime }(u^{n}(t))-g^{\prime}(u^{n-1}(t))\big)u_{x}^{n}(t)\big\Vert_{H^{s-1}}(e^{Ct}-1) \\
        &\leq & C \big\Vert \big(g^{\prime }(u^{n}(t))-g^{\prime}(u^{n-1}(t))\big)u_{x}^{n}(t)\big\Vert_{H^{s-1}}(e^{Ct}-1) \\
        &\leq & C \big\Vert u^{n}(t)-u^{n-1}(t)\big\Vert_{H^{s-1}}(e^{Ct}-1),
\end{eqnarray*}%
where we have used Lemma \ref{lem4.1}. Choosing $T_{1}<T_{0}$ now so that $e^{CT_{1}}-1\leq \frac{1}{2C}$ we see that for $\ t\leq T_{1},$
\begin{equation*}
    \Vert p^{n+1}(t)\Vert _{H^{s-1}}+\Vert q^{n+1}(t)\Vert_{H^{s-1}}
        \leq \frac{1}{2}\big(\Vert p^{n}(t)\Vert _{H^{s-1}}+\Vert q^{n}(t)\Vert _{H^{s-1}}\big)\leq \cdots \leq \frac{C}{2^{n}}.
\end{equation*}%
This shows that $(u^{n}, v^{n})$ forms a Cauchy sequence in $H^{s-1}$ and the limit $( u,v) $ will be a solution in $H^{s-1}.$ Finally,
as was done in the proof of Lemma \ref{lem3.1}, we obtain regularity, namely that $u,v\in C\big([ 0,T_{1}], H^{s}\big)$. Noting that the operator $K$ has a smoothing effect of order $r/2$, we get $u_{t}, v_{t}\in C\big([0,T_{1}],H^{s+\frac{r}{2}-1}\big)$ when $r<2$. This proves the
theorem below:
\begin{theorem}\label{theo4.2}
        Suppose $\beta$ satisfies (\ref{r}) for some $r\geq 0$.     Let $s>5/2$, $\sigma=\min \{ s,s+\frac{r}{2}-1\}$. Suppose $u_{0},v_{0}\in H^{s}$  with sufficiently small $\Vert u_{0}\Vert _{H^{s}}+\Vert v_{0}\Vert _{H^{s}}$. Then there exists  $T_{1}>0$ so that the Cauchy problem  (\ref{systema})-(\ref{systemb})  has a unique solution $\ u,v\in Y^{s, \sigma}=  C\big( [0, T_{1}] ,H^{s}\big) \cap C^{1}\big( [0,T_{1}],H^{\sigma}\big)$.
\end{theorem}
We conclude by stating the local well-posedness result for our original equation (\ref{nw}). Noting the relationship between the initial values for (\ref{nw}) and (\ref{systema})-(\ref{systemb}), we have:
\begin{corollary}\label{cor4.3}
        Suppose $\beta$ satisfies (\ref{r}) for some $r\geq 0$.     Let $s>5/2$, $\sigma=\min \{ s,s+\frac{r}{2}-1\}$. Suppose $u_{0},v_{0}\in H^{s}$  with sufficiently small $\Vert u_{0}\Vert _{H^{s}}+\Vert v_{0}\Vert _{H^{s}}$. Then there exists $T_{1}>0$ so that the Cauchy problem for (\ref{nw}), with initial data $u(x,0)=u_{0}(x), $ $u_{t}(x,0)=K( v_{0}(x))_{x}$ has a unique solution $u\in Y^{s, \sigma}=  C\big( [0, T_{1}] ,H^{s}\big) \cap C^{1}\big( [0,T_{1}],H^{\sigma}\big)$.
\end{corollary}
\begin{remark}\label{rem4.4}
    We finally note that when $K$ is invertible, namely $\widehat{\beta }(\xi)\geq c(1+\xi ^{2})^{-r/2},$ \ the condition $u_{t}(x,0)=K( v_{0}(x))_{x}$ amounts to $u_{t}(x,0)=( w_{1}(x))_{x}$ with arbitrary $w_{1}\in H^{\sigma}.$
\end{remark}

\setcounter{equation}{0}
\section{Blow up and existence time of solutions}\label{sec5}

In this section we will first determine conditions for finite-time blow-up of solutions to (\ref{systema})-(\ref{systemb})and then we will highlight the dependence of the life span of solutions upon the initial data. All the results in this section correspond to similar results for our original problem (\ref{nw}).

When $r\geq 2,$ it was shown in \cite{Duruk2010} that blow-up can only happen in the $L^{\infty }$ norm of $u$; namely, as long as $\Vert u(t)\Vert _{L^{\infty }}$ stays finite, there is global existence for any $s>1/2$. This is due to the fact that when $r\geq 2,$ (\ref{systema})-(\ref{systemb}) is an $H^{s}$-valued ODE system. In the general case $r\geq 0,$ the condition that $\Vert u(t)\Vert _{L^{\infty }}$ stays finite will not suffice as there may be also wave breaking. The phenomenon where $\Vert u(t)\Vert _{L^{\infty }}$ remains bounded but $\Vert u_x(t)\Vert _{L^{\infty }}$ and/or $\Vert u_{t}(t)\Vert _{L^{\infty }}$ become unbounded will be called as wave breaking.

Consider the solution $(u,v)$ of (\ref{systema})-(\ref{systemb}) with the initial data $u_{0},v_{0}$ in $H^{s}$, $s>5/2$.  Then, obviously, the pair $(u,v)$ satisfies (\ref{inhomA})-(\ref{inhomB}) with $w=g^{\prime }(u)$ and $F_{1}=F_{2}=G=0$.
From the estimate (\ref{e-com}) in the proof of Lemma \ref{lem2.2}  we have
\begin{equation*}
      \frac{d}{dt}\mathcal{E}_{s}^{2}
            \leq  C\Big(\Vert w_{t}\Vert_{L^{\infty}} \Vert u\Vert_{H^{s}}^{2}
                +\Vert w_{x}\Vert_{L^{\infty}} \Vert v\Vert_{H^{s}}\Vert u\Vert_{H^{s}}
                +\Vert v\Vert_{H^{s}} \big\Vert [ \Lambda^{s},w] u_{x}\big\Vert_{L^{2}}\Big).
\end{equation*}%
When $w=g^{\prime }(u)$, with the use of the stronger commutator estimate
\begin{equation*}
    \big\Vert [ \Lambda^{s},w] u_{x}\big\Vert _{L^{2}}
        \leq C\Big(\Vert w_{x}\Vert _{L^{\infty }}\Vert u_{x}\Vert_{H^{s-1}}+\Vert w\Vert _{H^{s}}\Vert u_{x}\Vert _{L^{\infty }}\Big)
\end{equation*}%
in Lemma \ref{lem2.1}, this gives
\begin{eqnarray}
    \frac{d}{dt}\mathcal{E}_{s}^{2}
        &\leq & C\Big(\Vert (g^{\prime }(u))_{t}\Vert_{L^{\infty }}\Vert u\Vert _{H^{s}}^{2}
        +\Vert (g^{\prime}(u))_{x}\Vert _{L^{\infty }}\Vert u\Vert _{H^{s}}\Vert v\Vert _{H^{s}} \nonumber \\
       && ~~ +\big(\Vert (g^{\prime}(u))_{x}\Vert _{L^{\infty }}\Vert u\Vert_{H^{s}}+\Vert g^{\prime }(u)\Vert _{H^{s}}\Vert u_{x}\Vert _{L^{\infty }}\big)\Vert v\Vert _{H^{s}}\Big). \label{masterEs}
\end{eqnarray}%
Now suppose that for some $T>0$,  $\Vert u(t)\Vert_{L^{\infty }}\leq  M$ and $\Vert u_{t}(t)\Vert_{L^{\infty }}+\Vert u_{x}(t)\Vert_{L^{\infty }}\leq  M_{1}$ for all $t\in [0, T]$. But with
\begin{equation*}
    C(M)=\sup_{\vert p \vert \leq M}\big(\vert g^{\prime}(p)\vert + \vert g^{\prime\prime}(p)\vert\big)
\end{equation*}
$\Vert g^{\prime }(u)\Vert _{L^{\infty }}\leq C(M)$ and $\Vert (g^{\prime }(u))_{t}\Vert _{L^{\infty}}+\Vert (g^{\prime }(u))_{x}\Vert_{L^{\infty }}\leq C(M)\big(\Vert u_{t}\Vert_{L^{\infty }}+\Vert u_{x}\Vert _{L^{\infty }}\big)\leq C(M)M_{1}$, and $\Vert g^{\prime }(u)\Vert _{H^{s}}\leq C(M)\Vert u\Vert_{H^{s}}$ for $t\in [0, T]$, where we have used Lemma \ref{lem4.1}. So we get
\begin{equation*}
    \frac{d}{dt}\mathcal{E}_{s}^{2}(t)\leq C(M)M_{1}\mathcal{E}_{s}^{2}(t),
\end{equation*}%
and Gronwall's inequality implies
\begin{equation*}
    \Vert u(t)\Vert^2_{H^{s}}+\Vert v(t)\Vert^2_{H^{s}}\approx \mathcal{E}_{s}^{2}(t)
        \leq C e^{C(M)M_{1}t}.
\end{equation*}%
This shows that whenever $\Vert u(t)\Vert _{L^{\infty }}$,  $\Vert u_{t}(t)\Vert _{L^{\infty }}$ and $\Vert u_{x}(t)\Vert _{L^{\infty }}$ stay finite, $\Vert u(t)\Vert^2_{H^{s}}+\Vert v(t)\Vert^2_{H^{s}}$ also stays finite; so to investigate blow-up versus global existence, one needs to control  the $L^{\infty }$  norms of $u$, $u_{x}$ and $u_{t}$. In other words;
\begin{lemma}\label{lem5.1}
    Let $s>5/2$. Then the solution of (\ref{systema})-(\ref{systemb}), with initial data in $H^s$, will blow up in finite time if and only if
    \begin{equation*}
        \Vert u(t)\Vert _{L^{\infty }}+\Vert u_{t}(t)\Vert_{L^{\infty }}+\Vert u_{x}(t)\Vert _{L^{\infty }}
    \end{equation*}%
    blows up in finite time.
\end{lemma}
Next we will consider the existence time of solutions to (\ref{systema})-(\ref{systemb}) in the case of a power-type nonlinearity,  $g(u)=\alpha u^{n+1}$ ($n \geq 1$); namely we show  how existence time depends on the initial data. Then, with  $w=g^{\prime }(u)=(n+1)\alpha u^{n}$ we have from (\ref{masterEs}):
\begin{equation}
    \frac{d}{dt}\mathcal{E}_{s}^{2}
        \leq  C\Big(\Vert \big(u^{n}\big)_{t}\Vert_{L^{\infty }}\Vert u\Vert _{H^{s}}^{2}
        +\Vert \big(u^{n}\big)_{x}\Vert _{L^{\infty }}\Vert u\Vert _{H^{s}}\Vert v\Vert _{H^{s}}+\Vert u\Vert^{n}_{H^{s}}\Vert u_{x}\Vert _{L^{\infty }}\Vert v\Vert _{H^{s}}\Big). \label{energyblow}
\end{equation}%
The first term in (\ref{energyblow}) is estimated as follows
\begin{equation*}
      \Vert \big(u^{n}\big)_{t}\Vert_{L^{\infty}} \leq C \Vert u\Vert^{n-1}_{L^{\infty}}\Vert u_{t}\Vert_{L^{\infty}}
                \leq C \Vert u\Vert^{n-1}_{H^{s-1}}\Vert u_{t}\Vert_{H^{s-1}}.
\end{equation*}%
But $\Vert u_{t}\Vert_{H^{s-1}}=\Vert Kv_{x}\Vert_{H^{s-1}}\leq C \Vert v_{x}\Vert_{H^{s-1}}\leq \Vert v\Vert_{H^{s}}$, so
\begin{equation*}
      \Vert \big(u^{n}\big)_{t}\Vert_{L^{\infty}} \leq C \Vert u\Vert^{n-1}_{H^{s}}\Vert v\Vert_{H^{s}}\leq C  \mathcal{E}_{s}^{n}.
\end{equation*}%
Similarly,
\begin{equation*}
      \Vert \big(u^{n}\big)_{x}\Vert_{L^{\infty}} \leq C \Vert \big(u^{n}\big)_{x}\Vert_{H^{s-1}}\leq C \Vert u\Vert^{n}_{H^{s}}\leq C  \mathcal{E}_{s}^{n}.
\end{equation*}%
Adding up,
\begin{equation*}
      \frac{d}{dt}\mathcal{E}_{s}^{2}(t)    \leq  C\mathcal{E}_{s}^{n+2}(t).
\end{equation*}%
This in turn implies that
\begin{equation*}
    \Vert u(t)\Vert_{H^{s}}+\Vert v(t)\Vert_{H^{s}}\approx \mathcal{E}_{s}(t)\leq \frac{\mathcal{E}_{s}(0)}{\big(1-Ct\mathcal{E}_{s}^{n}(0)\big)^{\frac{1}{n}}}.
\end{equation*}%
In other words, we have the following lemma establishing the relationship between the existence time and the initial data:
\begin{lemma}\label{lem5.2}
    The solution of (\ref{systema})-(\ref{systemb}) with $g(u)=\alpha u^{n+1}$ ($n \geq 1$)  and initial data $u_{0},v_{0}$ in $H^{s}$ exists at least for times
    \begin{equation*}
    t={\cal O}\Big(\frac{1}{\big(\Vert u_{0}\Vert _{H^{s}}+\Vert v_{0}\Vert _{H^{s}}\big)^{n}}\Big).
    \end{equation*}
\end{lemma}

We next investigate the dependence of the existence time on the size of the nonlinear term.  For this aim we consider the problem
\begin{eqnarray}
    && u^{\epsilon}_{t}=Kv^{\epsilon}_{x}~,\text{ \ \ \ \ \ \ \ \ \ \ \ \ \ \ \ \ \ \ \ \ \ }u^{\epsilon}(x,0)=u_{0}(x),~~\label{epsilona} \\
    && v^{\epsilon}_{t}=K\big(u^{\epsilon}+\epsilon^{n} (u^{\epsilon})^{n+1}\big)_{x},\text{ \ \ }v^{\epsilon}(x,0)=v_{0}(x),  \label{epsilonb}
\end{eqnarray}%
parametrized by a small positive parameter $\epsilon$  reflecting the strength of the nonlinearity. Using  again the estimate (\ref{e-com}) in Lemma \ref{lem2.2} with $F_{1}=F_{2}=G=0$ and keeping track of the extra $\epsilon^{n}$ in front of $(u^{\epsilon})^{n+1}$, we get
\begin{equation*}
    \Vert u^{\epsilon}(t)\Vert _{H^{s}}+\Vert v^{\epsilon}(t)\Vert_{H^{s}}\approx \mathcal{E}_{s}(t)
        \leq \frac{\mathcal{E}_{s}(0)}{\big(1-C t \epsilon^{n}\mathcal{E}^{n}_{s}(0)\big)^{\frac{1}{n}}}.
\end{equation*}%
That is,  the existence time of the time-dependent solution $(u,v)$ to (\ref{epsilona})-(\ref{epsilonb}) is at least of order ${\cal O}(1/\epsilon^{n})$ and to get the long-time existence of solutions, it suffices to take $\epsilon$ small enough.  So we have the following lemma establishing the long-time existence result:
\begin{lemma}\label{lem5.3}
    The solution  of (\ref{epsilona})-(\ref{epsilonb}) with $g(u)=\alpha u^{n+1}$ ($n \geq 1$)  and initial data $u_{0},v_{0}$ in $H^{s}$, exists at least for times
    \begin{equation*}
        t={\cal O}\Big(\frac{1}{\epsilon^{n}\big(\Vert u_{0}\Vert _{H^{s}}+\Vert v_{0}\Vert _{H^{s}}\big)^{n}}\Big).
    \end{equation*}
\end{lemma}

\setcounter{equation}{0}
\section{Convergence to the equation of classical elasticity}\label{sec6}

The object of this section is to study the relationship between the family of nonlocal wave equations and the   classical (local) elasticity model. We prove that as the kernel $\beta $ approaches the Dirac delta measure, not only (\ref{nw}) approaches formally the classical elasticity equation, (\ref{classical}), but also the corresponding solutions converge strongly to the solution of (\ref{classical}). The notion of this comparison result is in the sense of that in \cite{Erbay2020}.

We will consider the largest class of kernels, namely those satisfying (\ref{rbound}). Note that as $\widehat{\beta}$ is real, $\beta$ and $\widehat{\beta}$ are even functions.  We will further assume that
\begin{equation}
  \int_{\mathbb{R}} \beta dx=1~~~~\text{and} ~~~~ \int_{\mathbb{R}} x^{2}\vert \beta \vert dx<\infty.  \label{moment}
\end{equation}%
In the rest of the paper we focus on the power-type nonlinearity   $g(u)=u^{n+1}$ ($n \geq 1$). Using the transformation $(x, t, u) \rightarrow (x/\delta, t/\delta, \epsilon u^{\delta})$  with  small parameters $\delta >0$ and $\epsilon >0$ in (\ref{nw}), we now define the family of kernels as $\beta _{\delta }(x)=\frac{1}{\delta }\beta (\frac{x}{\delta })$  and the family of nonlocal nonlinear dispersive wave equations
\begin{equation}
    u^{\delta}_{tt}=\beta_{\delta }\ast \big(u^{\delta}+\epsilon^{n} (u^{\delta})^{n+1}\big)_{xx}.  \label{nwdelta}
\end{equation}%
Here, the parameter $\delta$ acts as a length scale in the problem and serves to measure the intensity of the dispersive effect. We note that, as $\delta \rightarrow 0$, the kernels $\beta_{\delta}$ will converge to the Dirac measure  in the distribution sense. Hence,  as the dispersion parameter $\delta$ tends to zero, (\ref{nwdelta}) will formally approach the classical elasticity equation
\begin{equation}
    u_{tt}=\big(u+\epsilon^{n} u^{n+1}\big)_{xx}  \label{claseps}
\end{equation}%
corresponding to the nondispersive case. Understanding this passage from (\ref{nwdelta}) to (\ref{claseps}) is  our main motivation here. We shall use $u^{\delta}$ and $u$ to denote the solutions to the initial value-problems defined for the nonlocal model (\ref{nwdelta}) and the local model  (\ref{claseps}), respectively. We will prove that  the difference between the corresponding  solutions of (\ref{nwdelta}) and (\ref{claseps})  remains small over a long time interval in a suitable Sobolev norm. That is, we will show that,  in the zero-dispersion limit, the solutions of  the classical elasticity equation are well approximated by the corresponding solutions of the nonlocal wave equations. Obviously, the solutions we consider are small-but-finite amplitude long wave solutions  as it is clear from the form of (\ref{nwdelta}).

First we note that the relationship between the Fourier transforms of $\beta _{\delta}$ and $\beta$ is   $\widehat{\beta _{\delta }}(\xi )=\widehat{\beta }(\delta \xi )$. So $\Vert \widehat{\beta _{\delta }}\Vert_{L^{\infty}}=\Vert \widehat{\beta }\Vert_{L^{\infty}}$ and thus the energy
estimates of the previous section will be uniform in $\delta$. Hence Lemma \ref{lem5.3} implies that the solution $u^{\delta }$ of the Cauchy problem for (\ref{nwdelta}) will all exist in a certain time interval $\ [ 0,\frac{T}{\epsilon^{n} }] $. Similarly,  when $\beta$ is the Dirac delta measure corresponding to (\ref{claseps}), we have existence of solution $u$  for times ${\cal O}(\frac{1}{\epsilon^{n} })$ hence we can assume that both $u$ and $u^{\delta}$ exist in the same time interval. Moreover, as $\sigma\geq s-1$ in Corollary \ref{cor4.3}, it is clear that both $u$ and $u^{\delta}$  belong to $C\big([0, \frac{T}{\epsilon^{n} }], H^{s}\big) \cap C^{1}\big([0, \frac{T}{\epsilon^{n}}], H^{s-1}\big)$. Our aim is to show that, for suitable initial data, $u^{\delta}$ converges to $u$ in $C\big([0,\frac{T}{\epsilon^{n} }], H^{s}\big) \cap C^{1}\big([0, \frac{T}{\epsilon^{n} }], H^{s-1}\big)$.

We convert (\ref{nwdelta}) to the system
\begin{eqnarray}
    u_{t}^{\delta}&=&K_{\delta }v^{\delta}_{x}, \label{sysdeltaA} \\
    v^{\delta}_{t}&=&K_{\delta}\big (u^{\delta}+\epsilon^{n} (u^{\delta})^{n+1}\big)_{x},\text{ }  \label{sysdeltaB}
\end{eqnarray}%
where $K_{\delta }v=\mathcal{F}^{-1}\big( k(\delta \xi )\widehat{v}(\xi)\big) $ while (\ref{classical}) corresponds to
\begin{eqnarray}
    u_{t}&=& v_{x},  \label{sysclasA} \\
    v_{t}&=& \big(1+(n+1)\epsilon^{n} u^{n}\big)u_{x}.  \label{sysclasB}
\end{eqnarray}%
We will use the estimates of Lemma \ref{lem2.2} to compare the solutions of these two systems.

We first derive an estimate on the operator $K_{\delta }-I$.
\begin{proposition}\label{prop6.1}
    For $v\in H^{s+\theta}$ and $0<\theta \leq 2$, we have
    \begin{equation}
        \Vert K_{\delta}v-v\Vert_{H^{s}}\leq C\delta ^{\theta}\Vert v\Vert_{H^{s+\theta }}. \label{K-estimate}
    \end{equation}
\end{proposition}
\begin{proof}
    The assumption (\ \ref{moment}) implies that $\widehat{\beta }^{\prime\prime }$ exists and is bounded. Also, since $\widehat{\beta }$ is real, $\beta $ and hence $\widehat{\beta }$ \ will be even. Clearly $k=(\widehat{\beta })^{1/2}$ will have the same properties. As $k( 0) =1$ and $k^{\prime }(0)=0,$ we have the Taylor estimate $\vert k(\xi)-1\vert \leq \frac{1}{2}\Vert k^{\prime \prime }\Vert _{L^{\infty }}\xi ^{2}$. Joining this estimate with $\vert k(\xi )-1\vert \leq \Vert k\Vert_{L^{\infty }}+1$, we have for any $0<\theta \leq 2$
    \begin{equation}
        \vert k(\xi )-1\vert
        = \vert k(\xi)-1\vert ^{1-\frac{\theta}{2} }\vert k(\xi )-1\vert^{\frac{\theta}{2} }
        \leq C\vert \xi \vert ^{\theta },  \label{theta}
    \end{equation}%
    where $C$ is a constant dependent on $\Vert k\Vert_{L^{\infty }}$ and $\Vert k^{\prime \prime }\Vert _{L^{\infty }}$.     Then
    \begin{eqnarray*}
        \Vert K_{\delta }v-v\Vert _{H^{s}}^{2}
            &=&\int_{\mathbb{R}} (k(\delta \xi )-1)^{2}( 1+\xi ^{2}) ^{s}\vert\widehat{v}(\xi )\vert ^{2}d\xi \\
            &\leq & C^{2}\delta^{2\theta }\int_{\mathbb{R}} \vert \xi \vert ^{2\theta}( 1+\xi ^{2})^{s}\vert \widehat{v}(\xi)\vert ^{2}d\xi
            \leq C^{2}\delta^{2\theta}\Vert v\Vert_{H^{s+\theta}}^{2}.
    \end{eqnarray*}
\end{proof}
\begin{theorem}\label{theo6.2}
    Suppose $\beta$ satisfies (\ref{rbound}) and (\ref{moment}). Let $u_{0},v_{0}\in H^{s+\theta}$, $s>5/2$, $0<\theta \leq 2$. Let $u^{\delta},v^{\delta}, u,v\in C\big( [ 0,\frac{T}{\epsilon^{n}}], H^{s+\theta}\big) \cap C^{1}\big([0,\frac{T}{\epsilon^{n}}], H^{^{s+\theta-1 }}\big) $ satisfy (\ref{sysdeltaA})-(\ref{sysdeltaB}) and (\ref{sysclasA})-(\ref{sysclasB}) with the  initial data $u^{\delta}(x,0)=u(x,0)=u_{0}(x)$ and $v^{\delta}(x,0)=v(x,0)=v_{0}(x)$. Then we have, for all $t\in [ 0,\frac{T}{\epsilon^{n} }] $,
        \begin{equation}
            \Vert u^{\delta}(t)-u(t)\Vert _{H^{s-1}}+\Vert v^{\delta}(t)-v(t)\Vert _{H^{s-1}}\leq C \delta^{\theta} t e^{C\epsilon^{n}t},
             \label{maintheo}
        \end{equation}
    where $C$ is a constant independent of $\delta$ and $\epsilon$.
\end{theorem}
\begin{proof}
    By Lemma \ref{lem2.2}, we have  uniform bounds on $\Vert u^{\delta}\Vert _{X^{s+\theta}}$, $\Vert v^{\delta}\Vert _{X^{s+\theta}}$, $\Vert u\Vert_{X^{s+\theta}}$ and $\Vert v\Vert _{X^{s+\theta}}$. Let $p=u^{\delta}-u$, $q=v^{\delta}-v$. Then the diferences $p$ and $q$ solve the system
    \begin{eqnarray*}
        p_{t} &=&K_{\delta }v_{x}^{\delta}-v_{x} \\
        q_{t} &=&K_{\delta }u_{x}^{\delta}+(n+1)\epsilon^{n} K_{\delta }\big((u^{\delta})^{n}u_{x}^{\delta}\big)
                    -u_{x}-(n+1)\epsilon^{n} u^{n}u_{x},
    \end{eqnarray*}%
    or equivalently the system
    \begin{eqnarray*}
        p_{t} &=&K_{\delta}q_{x}+K_{\delta}v_{x}-v_{x} \\
        q_{t} &=&K_{\delta}p_{x}+K_{\delta}u_{x}-u_{x}+(n+1)\epsilon^{n} \Big( K_{\delta}\big((u^{\delta})^{n}p_{x}\big)
                    +  K_{\delta}\big((u^{\delta})^{n}u_{x}\big)                      \\
                    && -(u^{\delta})^{n}u_{x} + (u^{\delta})^{n}u_{x}                    -u^{n}u_{x}\Big).
    \end{eqnarray*}%
    We have the initial-value problem
    \begin{eqnarray}
        p_{t} &=&K_{\delta}q_{x}+F_{1}, \text{\ \ \ \ \ \ \ \ \ \ \ \ \ \ \ \ \ \ \ \ \ \ \ \ \ \ \ \ \ \ \ \ \ \ \ \ \ \ \ \ \ } p(x,0)=0,  \label{elasA}\\
        q_{t} &=&K_{\delta}\Big(\big(1+(n+1)\epsilon^{n} (u^{\delta})^{n}\big)p_{x}\Big)+F_{2}+G, \text{\ \ \ \ } q(x,0)=0,  \label{elasB}
    \end{eqnarray}%
    where the residual terms $F_{1}$, $F_{2}$ and $G$ are given by
    \begin{eqnarray*}
        F_{1} &=&K_{\delta}v_{x}-v_{x} \\
        F_{2} &=& K_{\delta}u_{x}-u_{x}+(n+1)\epsilon^{n} \Big(K_{\delta}\big((u^{\delta})^{n}u_{x}\big)
                -(u^{\delta})^{n}u_{x}\Big) \\
        G &=& (n+1)\epsilon^{n} \big( (u^{\delta})^{n}-u^{n}\big)u_{x}.
    \end{eqnarray*}%
    To complete the proof we need to  estimate the residual terms. By Proposition \ref{prop6.1}, we have the following bounds on $F_{1}$, $F_{2}$ and $G$:
    \begin{eqnarray*}
        \Vert F_{1}(t)\Vert _{H^{s-1}}
            &=& \big\Vert K_{\delta}v_{x}(t)-v_{x}(t)\big\Vert _{H^{s-1}} \\
            &\leq & C\delta^{\theta }\Vert v_{x}(t)\Vert _{H^{s+\theta-1 }}
                    \leq  C\delta^{\theta }\Vert v(t)\Vert _{H^{s+\theta}}\leq C\delta^{\theta} \\
        \Vert F_{2}(t)\Vert _{H^{s-1}}
            &\leq & \big\Vert K_{\delta}u_{x}(t)-u_{x}(t)\big\Vert _{H^{s-1}}
                +(n+1)\epsilon^{n} \big\Vert K_{\delta}((u^{\delta})^{n}u_{x})(t)
                        -((u^{\delta})^{n}u_{x})(t)\big\Vert _{H^{s-1}}      \\
            &\leq& C\delta^{\theta}+C\epsilon^{n} \delta^{\theta}\Vert \big((u^{\delta})^{n}u_{x}\big)(t)\Vert_{H^{s+\theta-1}} \\
            &\leq& C\delta^{\theta}+C\epsilon^{n} \delta^{\theta}\Vert u(t)\Vert_{H^{s+\theta}}
                    \leq C\delta ^{\theta} \\
    \Vert G(t)\Vert _{H^{s-1}}
        &=&(n+1)\epsilon^{n} \big\Vert \big((u^{\delta}(t))^{n}-(u(t))^{n}\big)u_{x}(t)\big\Vert_{H^{s-1}} \\
        &\leq& C\epsilon^{n} \big\Vert (u^{\delta}(t))^{n}-(u(t))^{n}\big\Vert_{H^{s-1}}\Vert u_{x}(t)\Vert_{H^{s-1}} \\
        &\leq& C\epsilon^{n} \big\Vert (u^{\delta}(t))^{n}-(u(t))^{n}\big\Vert_{H^{s-1}}\Vert u(t)\Vert_{H^{s}}\\
        &\leq& C\epsilon^{n} \big\Vert u^{\delta}(t)-u(t)\big\Vert_{H^{s-1}}
        =C\epsilon^{n}\Vert p(t)\Vert_{H^{s-1}}
    \end{eqnarray*}
    for some constant $C$ independent of $\delta$. Combining these estimates with the one in (\ref{e-com}) we get
     \begin{equation}
        \frac{d}{dt}\mathcal{E}_{s-1}^{2}(t)
            \leq C\epsilon^{n}\mathcal{E}_{s-1}^{2}(t)+C \delta^{\theta}\mathcal{E}_{s-1}(t). \label{elasEN}
    \end{equation}%
     Applying Gronwall's inequality yields (\ref{maintheo}). This completes the proof.
\end{proof}

The above theorem is the main result of our work. It says that the error in the approximation of solutions to  (\ref{claseps}) by the corresponding solutions of (\ref{nwdelta}) will be small as $\delta$ becomes small enough with some convergence rate $\theta \in (0, 2]$ and it  is controlled on the time scale ${\cal O}(1/\epsilon^{n})$. Also, the convergence becomes slower for larger times.

\setcounter{equation}{0}
\section{Convergence from the discrete lattice model to the continuous model}\label{sec7}

As a continuation of the considerations in the previous section,  the problem of  searching the continuous limit of the discrete lattice dynamic model is interesting in its own right. So, our objective in this section is to  provide a rigorous justification of the convergence of the FPUT equation (which is a member of the class  (\ref{nw}) for  the triangular kernel to be described below)  to the classical elasticity equation in the limit of vanishing nonlocality. In other words, we now show that, as the lattice spacing approaches zero, the FPUT discrete lattice model and the classical elasticity model behave qualitatively similarly, provided the initial data is properly imposed.  As it is expected, in the long-wave limit the discreteness does not prevail and the predictions of the discrete and  continuum models agree to within an order of neglected terms.

To this end we consider  the triangular kernel  defined by
\begin{equation*}
    \beta ^{T}(x)=\left\{
        \begin{array}{cc}
        1-|x| & \text{ for }~|x|\leq 1, \\
        0 & \text{ for }~|x|>1,
        \end{array}%
        \right.
\end{equation*}%
for which  $ \widehat{\beta^{T}}(\xi )=\frac{4}{\xi ^{2}}\sin ^{2}( \frac{\xi }{2})$. Recall that $k^{T}(\xi)=\big( \widehat{\beta^{T}}(\xi)\big)^{1/2}$.  Then (\ref{nwdelta}) becomes the differential-difference equation
\begin{equation}
    u^{\delta}_{tt}=\Delta _{\delta }^{d}\big(u^{\delta}+\epsilon^{n}(u^{\delta})^{n+1}\big)  \label{lattice}
\end{equation}%
if we use  the difference operator
\begin{equation}
    \Delta _{\delta }^{d}g(u)=\frac{1}{\delta^{2}}\Big(g\big(u(x+\delta,t)\big)-2g\big(u(x,t)\big)+g\big(u(x-\delta,t)\big)\Big).     \label{lap}
\end{equation}
Under the transformation $u^{\delta}(x,t)=\big(w^{\delta}(x+\delta,t)-w^{\delta}(x,t)\big)/\delta$, the differential-difference equation (\ref{lattice}) becomes the famous FPUT equation
\begin{equation}
    w^{\delta}_{tt}(x,t)=\Delta_{\delta }^{d}w^{\delta}(x,t)+\frac{\epsilon^{n}}{\delta^{n+2}} \Big(\big(w^{\delta}(x+\delta ,t)-w^{\delta}(x,t)\big)^{n+1}-\big(w^{\delta}(x,t)-w^{\delta}(x-\delta,t)\big)^{n+1}\Big)  \label{fput}
\end{equation}%
The discrete lattice model (\ref{lattice}), or equivalently  (\ref{fput}), describes longitudinal vibrations of a one-dimensional doubly-infinite chain  of identical particles with unit mass. The lattice consists of particles equally spaced at distance $\delta$  and it takes into account only nearest neighbor interactions  (this is the celebrated FPUT lattice \cite{Fermi1955} with the potential function $V(u)=\frac{1}{2}\epsilon^{2}u^{2}+\frac{\epsilon^{n+2}u^{n+2}}{n+2}$).  In  the above equations,  $w^{\delta}(x,t)$  stands for the displacement of the particle at $x$ from its equilibrium position at time $t$ while $u^{\delta}(x,t)$ corresponds to the strain (the relative displacement  between locations of two adjacent particles at $x+\delta$ and $x$).

We are interested in investigating asymptotic behavior of (\ref{lattice}) as $\delta \rightarrow 0$, which allows us to pass from the discrete description to the continuous description. If we take the initial velocity as $u^{\delta}_{t}(x,0)=K_{\delta}^{T}(v_{0}(x))_{x}$ for some suitable $v_{0}$, the long-time existence result  of Section \ref{sec5} (Lemma \ref{lem5.3}) will hold. Through the Fourier transform one gets $K^{T}u=\chi_{(-\frac{1}{2},\frac{1}{2}) }\ast u$ where $\chi_{(-\frac{1}{2},\frac{1}{2}) }$ denotes the characteristic function of the interval $(-\frac{1}{2},\frac{1}{2})$; $K_{\delta }^{T}$ is the convolution with $\frac{1}{\delta }\chi_{(-\frac{\delta}{2},\frac{\delta}{2}) }$. Then $K_{\delta }^{T}D_{x}$ is the difference operator. So, for the initial values $u^{\delta}(x,0)=u_{0}(x)$ and
\begin{equation}
   u^{\delta}_{t}(x,0)=K_{\delta }^{T}\big(v_{0}(x)\big)_{x}
                                =\frac{1}{\delta }\big(v_{0}(x+\frac{\delta }{2})-v_{0} (x-\frac{\delta }{2})\big),     \label{referee}
\end{equation}%
(\ref{lattice}) has a solution over long times. \ As would be expected, $K_{\delta }^{T}(v_{0}(x))_{x}$ is just the discretized derivative of $v_{0}$ and $\lim_{\delta \rightarrow 0}K_{\delta }^{T}(v_{0}(x))_{x}=v_{0}^{\prime}(x).$ In fact from the Taylor theorem
\begin{eqnarray*}
    K_{\delta }^{T}(v_{0}(x))_{x}-v_{0}^{\prime }(x)
        &=&\frac{1}{2\delta }\int_{-\frac{\delta }{2}}^{\frac{\delta }{2}}(\frac{\delta }{2}-\vert y\vert )^{2}v^{\prime\prime\prime}_{0}(x+y)dy, \\
        \big\vert K_{\delta }^{T}(v_{0}(x))_{x}-v_{0}^{\prime }(x)\big\vert^{2}
        &\leq &\frac{1}{4\delta ^{2}} \Big(\int_{-\frac{\delta }{2}}^{\frac{\delta }{2}}(\frac{\delta }{2}-\vert y\vert )^{2}v^{\prime\prime\prime}_{0}(x+y)dy\Big)^{2} .
\end{eqnarray*}%
Then
\begin{equation*}
    \big\Vert K_{\delta }^{T}(v_{0})_{x}-v_{0}^{\prime}\big\Vert_{L^{2}}^{2}
            \leq \frac{1}{4\delta ^{2}}\Vert v_{0}^{\prime \prime\prime}\Vert_{L^{2}}^{2}
                \Big( \int_{-\frac{\delta }{2}}^{\frac{\delta }{2}}(\frac{\delta }{2}-\vert y\vert )^{2}dy\Big)^{2}
            \leq C\delta^{4}\Vert v_{0} \Vert _{H^{3}}^{2}.
\end{equation*}%
Since $K_{\delta }^{T}$ commutes with the operator $\Lambda ^{s}=(1-D_{x}^{2}) ^{s},$ we have  for any $s\geq 1$,
\begin{equation*}
    \big\Vert K_{\delta }^{T}(v_{0})_{x}-v_{0}^{\prime}\big\Vert_{H^{s-1}}
        =\big\Vert K_{\delta }^{T}\Lambda^{s-1}(v_{0})_{x}-\Lambda^{s-1}v_{0}^{\prime}\big\Vert _{L^{2}}
        \leq C\delta ^{2}\Vert \Lambda^{s-1}v_{0} \Vert_{H^{3}}
        \leq C\delta ^{2}\Vert v_{0} \Vert_{H^{s+2}}.
\end{equation*}%
 We now compare the solution $u^{\delta}$  of (\ref{lattice}) with initial data
\begin{equation}
    u^{\delta}(x,0)=u_{0}(x), \qquad u^{\delta} _{t}(x,0)
                =\frac{1}{\delta }\big( v_{0}(x+\frac{\delta }{2})-v_{0} (x-\frac{\delta }{2})\big) , \label{latticeinitial}
\end{equation}%
to the solution $u$ of the classical (local) elasticity equation, (\ref{claseps}), with the initial data
\begin{equation}
    u(x,0)=u_{0}(x),\quad u_{t}(x,0)=(v_{0}(x))_{x}, \label{classicalinitial}
\end{equation}%
where $u_{0}, v_{0} \in H^{s+2}$.   We first define the differences $p=u^{\delta}-u$ and $q=u^{\delta}_{t}-u_{t}$. Then, as in the previous section, $p, q$  satisfy the system (\ref{elasA})-(\ref{elasB}) with the initial data
\begin{equation*}
    p(x,0)=0,\quad q(x,0)=\frac{1}{\delta}\big(v_{0} (x+\frac{\delta }{2})-v_{0} (x-\frac{\delta }{2})\big)-(v_{0}(x))_{x}.
\end{equation*}%
However, Theorem \ref{theo6.2} will not directly apply due to nonzero initial data related to $q$. So, using Gronwalls' inequality and the fact that $\mathcal{E}_{s-1}(0)\approx\Vert q(0)\Vert_{H^{s-1}}\leq C\delta^{2}$, we conclude from (\ref{elasEN}) that
\begin{equation*}
    \big\Vert u^{\delta}(t) -u(t) \big\Vert _{H^{s-1}}+\big\Vert u^{\delta}_{t}(t) -u_{t}(t) \big\Vert _{H^{s-1}}\leq C\delta ^{2}(1+t)\text{ \ \ for all \ } t\leq \frac{T}{\epsilon^{n}}.\text{ }
\end{equation*}
With the above preparatory work, we can now state the following comparison result:
\begin{theorem}\label{theo7.1}
    Suppose $u_{0},v_{0}\in H^{s+2}$, $s>5/2$. Let $u^{\delta},v^{\delta},u,v\in C\big( [ 0,\frac{T}{\epsilon^n}], H^{s+2}\big) \cap C^{1}\big([0,\frac{T}{\epsilon^n}], H^{^{s+1}}\big) $ satisfy (\ref{lattice}) and (\ref{claseps}) with the initial data given by (\ref{latticeinitial}) and (\ref{classicalinitial}). Then we have, for all $t\in [ 0,\frac{T}{\epsilon^n}] $
        \begin{equation*}
             \big\Vert u^{\delta}(t) -u(t) \big\Vert_{H^{s-1}}+\big\Vert u^{\delta}_{t}(t) -u_{t}(t) \big\Vert _{H^{s-1}}\leq C\delta ^{2}(1+t)e^{C\epsilon^{n}t},
        \end{equation*}
    where $C$ is a constant independent of $\delta$ and $\epsilon$.
\end{theorem}
\begin{remark}
It is worth mentioning that the proper imposition of the initial conditions is crucial in establishing the above comparison result. The comparison result in Theorem \ref{theo7.1} is based on  the initial conditions   (\ref{latticeinitial}) and (\ref{classicalinitial})  given for $u^{\delta}$ and $u$, respectively. We note from (\ref{latticeinitial}) and (\ref{classicalinitial})  that while $u^{\delta}$ and $u$ are initially equal to the same function, the initial values of their derivatives $u^{\delta}_{t}$ and $u_{t}$  are different.  This is because the initial condition  stated in  Corollary \ref{cor4.3} for  $u^{\delta}_{t}$ involves  the pseudo-differential operator  $K$. So   $u^{\delta}_{t}(x,0)$ must be written properly in order to apply the existence result  given in Corollary \ref{cor4.3}  for the second-order equation  (\ref{nw}).
\end{remark}

\bibliographystyle{plainnat}
\bibliography{references}

\end{document}